\documentclass[10pt]{amsart}
\usepackage{amsmath,amssymb,amsthm,amscd,mathrsfs}
\usepackage[all]{xy}
\usepackage{color} 
\CompileMatrices
\numberwithin{equation}{section}
\newtheorem{prop}{Proposition}[section]
\newtheorem{theo}[prop]{Theorem}
\newtheorem{lemm}[prop]{Lemma}
\newtheorem{coro}[prop]{Corollary}

\numberwithin{equation}{section}

\newcommand{\be}{\begin{equation}}
\newcommand{\ee}{\end{equation}}

\newcommand\IZ{\mathbb {Z}}

\newcommand\IQ{\mathbb {Q}}
\newcommand{\IC}{\mathbb{C}}
\newcommand{\IE}{\mathbb{E}}
\newcommand{\IR}{\mathbb{R}}

\newcommand{\ba}{\begin{array}}
\newcommand{\ea}{\end{array}}

\newcommand{\CX}{{\mathcal X}}

\newcommand{\IH}{{\mathbb H}} 
\newcommand{\bal}{\begin{aligned}}
\newcommand{\eal}{\end{aligned}}

\newcommand{\mfd}{{\mathfrak d}}

\newcommand{\mfg}{{\mathfrak g}}

\newcommand{\mfh}{{\mathfrak h}}

\newcommand{\mfm}{{\mathfrak{M}}}

\newcommand{\longto}{\longrightarrow}

\newcommand{\CO}{{\mathcal O}}

\newcommand{\CE}{{\mathcal E}}
\newcommand{\CA}{{\mathcal A}}

\newcommand{\CC}{{\mathcal C}}

\newcommand{\obj}{{\mathfrak {Ob}}}

\newcommand{\mfcf}{{\mathfrak F}}

\newcommand{\mfcx}{{\mathfrak X}}

\title{Chamber Structure and Wallcrossing in the ADHM Theory of Curves II}
\author{Wu-yen Chuang, Duiliu-Emanuel Diaconescu, Guang Pan}

\begin{document}

\begin{abstract}
This is the second part of a project concerning variation of stability and 
chamber structure for ADHM invariants of curves. Wallcrossing formulas 
for such invariants are derived using the theory of stack function 
Ringel-Hall algebras constructed by Joyce and 
 the theory of generalized Donaldson-Thomas 
invariants of Joyce and Song. Some applications are presented, including 
strong rationality for local stable pair invariants of higher genus curves, and 
comparison with wallcrossing formulas of Kontsevich and Soibelman, and 
the halo formula of Denef and Moore. 
\end{abstract} 
\maketitle 

\tableofcontents

\section{Introduction}\label{intro}
Let $X$ be a smooth projective curve over $\IC$, $\CO_X(1)$ a very ample line bundle 
on $X$, and $M_1,M_2$ two line bundles 
on $X$ so that $M_1\otimes M_2\simeq K_X^{-1}$.
An ADHM sheaf $\CE$ on $X$ with twisting data 
$(M_1,M_2)$ is a coherent $\CO_X$-module 
$E$ decorated by morphisms 
\[
\Phi_{i}:E\otimes_X M_i \to E, \qquad 
\phi:E\otimes_X M_1\otimes_X M_2 \to \CO_X, \qquad 
\psi:\CO_X \to E
\]
with $i=1,2$, satisfying the ADHM relation 
\be\label{eq:ADHMrelation}
\Phi_1\circ(\Phi_2\otimes 1_{M_1}) - \Phi_2\circ (\Phi_1\otimes 1_{M_2}) 
+ \psi\circ \phi =0.
\ee
An  ADHM sheaf $\CE$ 
will be said to be of type $(r,e)\in \IZ_{\geq 0}\times \IZ$ if $E$ has 
rank $r\in \IZ_{\geq 0}$ and degree $e\in \IZ$. 

A triple $(E,\Phi_1,\Phi_2)$ 
with $\Phi_1,\Phi_2$ morphisms of $\CO_X$-modules as above satisfying 
relation \eqref{eq:ADHMrelation} for $\phi=0$, $\psi=0$, 
 will be called a Higgs sheaf on $X$ with coefficient 
sheaf $M_1\oplus M_2$.  

The following construction results concerning moduli spaces of 
ADHM sheaves were proved in the first part of this work \cite{chamberI}.
\begin{itemize}
\item There exists a stability condition for ADHM sheaves depending 
on a real parameter $\delta\in \IR$ \cite[Def. 2.1]{chamberI}, \cite[Def. 2.2]{chamberI}
so that for fixed $(r,e)\in \IZ_{\geq 1}\times \IZ$ the are finitely many 
critical stability parameters dividing the real axis into chambers. 
The set of $\delta$-semistable ADHM sheaves is constant within each chamber, 
and strictly semistable objects may exist only if $\delta$ takes a critical value. 
The origin $\delta=0$ is a critical value for all $(r,e)\in \IZ_{\geq 1}\times \IZ$.
\item For fixed $(r,e)\in \IZ_{\geq 1}\times \IZ$ 
and $\delta\in \IR$ there is an 
algebraic moduli stack of finite type over $\IC$
$\mfm_\delta^{ss}(\CX,r,e)$ of $\delta$-semistable 
locally free ADHM sheaves. If $\delta\in \IR$ is noncritical, 
$\mfm_\delta^{ss}(\CX,r,e)$ is  a quasi-projectve scheme 
equipped with a perfect obstruction theory \cite[Thm 1.2]{chamberI}, 
\cite[Thm 1.4]{chamberI}.
\item For fixed $(r,e)\in \IZ_{\geq 1}\times \IZ$ 
and $\delta\in \IR$ there is a natural algebraic torus ${\bf S}=\IC^\times$ action on 
the moduli stack $\mfm_\delta^{ss}(\CX,r,e)$ which acts on $\IC$-valued 
points by scaling the morphisms $(\Phi_1,\Phi_2)\to (t^{-1}\Phi_1,t\Phi_2)$, 
$t\in {\bf S}$. If $\delta$ is noncritical \cite[Thm 1.5]{chamberI} 
proves that the stack theoretic fixed locus $\mfm_\delta^{ss}(\CX,r,e)^{\bf S}$
is proper over $\IC$. Therefore residual ADHM 
invariants $A_\delta^{\bf S}(r,e)$
are defined by equivariant virtual integration in each stability 
chamber \cite[Def. 1.8]{chamberI}. 
\item For  $(r,e)\in \IZ_{\geq 1}\times \IZ$ there exists a critical value 
$\delta_M\in \IR_{>0}$ so that for any $\delta>\delta_M$,
$\mfm_\delta^{ss}(\CX,r,e)$
is isomorphic to the moduli space of stable pairs of Pandharipande and 
Thomas \cite{stabpairs-I} 
on the total space of the rank two bundle $M_1^{-1} \oplus M_2^{-1}$ 
on $X$. This identification includes the equivariant perfect obstruction theories 
establishing an equivalence between local stable pair theory and 
asymptotic ADHM theory (see \cite[Thm. 1.11]{modADHM} and \cite[Cor. 1.12]{modADHM}
for precise statements.)
\end{itemize}

The present paper represents the second part of this work. Its main goal is 
to derive wallcrossing 
formulas for the ADHM invariants $A_\delta^{\bf S}(r,e)$ using Joyce's stack function algebra theory and the theory of generalized Donaldson-Thomas invariants of Joyce and Song.
Moreover, it will be also shown that these formulas  
imply the BPS rationality conjecture formulated by Pandharipande and Thomas 
in \cite{stabpairs-I} for local stable pair invariants of curves. 

Similar results have been obtained in \cite{generating} for stable 
pair invariants of smooth projective Calabi-Yau threefolds defined via the 
the stack theoretic topological Euler character introduced by Joyce in 
\cite{J-motivic, J-II}. Moreover the wallcrossing formula relating stable pair and 
Donaldson-Thomas theory has been derived for the same type of invariants 
in \cite{curvecount,Hilbpairs}. The moduli spaces involved in the local 
construction considered admit natural Chern-Simons functionals, 
as explained in \cite[Sect. 7]{chamberI},
making the theory  
of Joyce and Song applicable to virtual residual stable pair invariants. 
Note also that analogous results have been subsequently obtained in 
\cite{walltwisted} for more general twisted quiver bundles on curves. 

\subsection{Wallcrossing formula}
Let $\delta_c \in \IR_{\geq 0}$ be a critical stability parameter of type 
 $(r,e)\in \IZ_{\geq 1}\times \IZ$, possibly zero, 
 and $\delta_+ >\delta_c$, $\delta_-<\delta_c$ be 
 stability parameters so that there are no critical stability parameters of type 
 $(r,e)$ in the 
 interval $[\delta_-, \delta_+]$.
 In order to simplify the formulas, we will denote the numerical 
 invariants by ${\alpha}=(r,e)$, and use the notation
 \[
 \mu_{\delta}({\alpha}) = {e+\delta \over r}, \qquad 
 \mu({\alpha})= {e\over r}
 \]
 for any ${\alpha}=(r,e)$ with $r\geq 1$, and any $\delta\in \IR$. 
 
 For fixed ${\alpha}=(r,e)$, $\delta_c\geq 0$ and $l\in \IZ_{\geq 2}$ let 
  ${\sf S}_{\delta_c}^{(l)}({\alpha})$
 be the set of all ordered decompositions  
 \be\label{eq:decompcondA}
 {\alpha}={\alpha}_1+\cdots +{\alpha}_{l}, \qquad {\alpha}_i=(r_i,e_i)\in 
 \IZ_{\geq 1}\times \IZ, \qquad 
 i=1,\ldots, l
 \ee
 satisfying
 \be\label{eq:decompcondB}
 \mu({\alpha}_1)=\cdots =\mu({\alpha}_{l-1}) =\mu_{\delta_c}({\alpha}_l)=\mu_{\delta_c}({\alpha}).
 \ee
  Note that the union ${\sf S}_{\delta_c}({\alpha})=\bigcup_{l\geq 2} 
 {\sf S}_{\delta_c}^{(l)}({\alpha})$ 
 is a finite set for fixed $\delta_c\geq 0$. 
  Then the following theorem is proven in section (\ref{wallsubsection}). 
  \begin{theo}\label{wallcrossingthmB}
$(i)$ The following wallcrossing formula holds for $\delta_c>0$ 
\be\label{eq:wallformulaAA}
\begin{aligned}
& A^{\bf S}_{\delta_+}({\alpha})-A^{\bf S}_{\delta_-}({\alpha}) = \\
& \mathop{\sum_{l\geq 2}}_{} {(-1)^{l-1}\over (l-1)!}\
\mathop{\sum_{({\alpha}_1, \ldots, {\alpha}_l)\in {\sf S}_{\delta_c}^{(l)}({\alpha})}}_{}
A^{\bf S}_{\delta_-}({\alpha}_l) 
\prod_{j=1}^{l-1} 
[(-1)^{e_{j}-r_{j}(g-1)}
(e_{j}-r_{j}(g-1)) H({\alpha}_j)].\\
\end{aligned}
\ee
$(ii)$ The following wallcrossing formula holds for $\delta_c=0$. 
\be\label{eq:wallformulaAB}
\begin{aligned}
& A^{\bf S}_{\delta_+}({\alpha})-A^{\bf S}_{\delta_-}({\alpha}) = \\
& \mathop{\sum_{l\geq 2}}_{} {(-1)^{l-1}\over (l-1)!}\
\mathop{\sum_{({\alpha}_1, \ldots, {\alpha}_l)\in {\sf S}_{0}^{(l)}({\alpha})}}_{}
A^{\bf S}_{\delta_-}({\alpha}_l) \prod_{j=1}^{l-1} 
[(-1)^{e_{j}-r_{j}(g-1)}
(e_{j}-r_{j}(g-1))H({\alpha}_j)]\\
& + \mathop{\sum_{l\geq 1}}_{} {(-1)^l\over l!}\
\mathop{\sum_{({\alpha}_1, \ldots, {\alpha}_l)\in {\sf S}_{0}^{(l)}({\alpha})}}_{}
  \prod_{j=1}^{l} 
[(-1)^{e_{j}-r_{j}(g-1)}
(e_{j}-r_{j}(g-1)) H({\alpha}_j)]\\
\end{aligned}
\ee
Moreover, if $g\geq 1$, the right hand sides of equations \eqref{eq:wallformulaAA},
\eqref{eq:wallformulaAB} vanish. 
\end{theo}
Here $H({\alpha})$ are  generalized Donaldson-Thomas type 
invariants for Higgs sheaves with numerical invariants ${\alpha}=(r,e)$ on $X$  
defined in section (\ref{wallsubsection}).

Some applications of Theorem (\ref{wallcrossingthmB}) are presented below. 
For any $(r,e)\in \IZ_{\geq 1}\times\IZ$ let $A_{\pm\infty}^{\bf S}(r,e)$ denote the ADHM invariants 
in the asymptotic chambers $\delta>>0$, $\delta<<0$ respectively. Then let 
\be\label{eq:inftygenfct}
Z_{\pm\infty}(q)_r = \mathop{\sum_{e\in \IZ}}_{} q^{e-r(g-1)} A_{\pm\infty}^{\bf S}(r,e). 
\ee
be the formal generating function of such invariants for fixed rank $r\geq 1$.
 According to \cite[Cor. 1.12]{modADHM}
$Z_{+\infty}(q)_r$ is the generating function of degree $r$ local stable pair invariants 
of the data $\CX=(X,M_1,M_2)$. Note that \cite[Lemm. 2.3]{chamberI} implies that 
$A_{+\infty}^{\bf S}(r,e)=A_{-\infty}^{\bf S}(r,-e+2r(g-1))$ for all $(r,e)\in \IZ_{\geq 1}\times\IZ$.
On the other hand, for curves $X$ of genus $g\geq 1$, Theorem (\ref{wallcrossingthmB}) 
implies that $A_{+\infty}^{\bf S}(r,e)=A_{-\infty}^{\bf S}(r,e)$ for all $(r,e)\in \IZ_{\geq 1}\times\IZ$
since in this case the invariants $H(r,e)$ are zero.
Therefore the following holds.
\begin{coro}\label{wallcrossingthmC}
If $g\geq 1$, $Z_{+\infty}(q)_r$
is a polynomial in $q, q^{-1}$ invariant under $q\leftrightarrow q^{-1}$.
\end{coro}
This implies that the strong rationality conjecture formulated in \cite{stabpairs-I} 
holds for the local stable pair theory of curves of genus $g\geq 1$. 

Further applications of Theorem (\ref{wallcrossingthmB}) are presented in 
section (\ref{KSsect}) where it is shown that the wallcrossing formula \eqref{eq:wallformulaAA} 
is in agreement with the wallcrossing formula of Kontsevich and Soibelman 
\cite{wallcrossing}. Moreover, it is also shown that formula \eqref{eq:wallformulaAA} 
is in agreement with the halo wallcrossing formula for D6-D2-D0 bound states 
derived by Denef and Moore  \cite[Sect. 6.1.2]{DM-split} 
using supergravity arguments. A more elaborate application of the above wallcrossing formulas 
to the cohomology of the moduli spaces of Hitchin pairs on curves is presented in 
\cite{wall-pairs}.

{\it Acknowledgements}. 
D.-E. D. would like to thank Arend Bayer,
Ugo Bruzzo, Daniel Jafferis, Dominic Joyce, Jan Manschot, Greg Moore, Sven Meinhart,
Kentaro Nagao, Alexander Schmitt, Andrei Teleman, Yukinobu Toda, and especially 
Ron Donagi, Liviu Nicolaescu and Tony Pantev for very helpful discussions, and 
correspondence and to 
Ionut Ciocan-Fontanine, Bumsig Kim and Davesh Maulik for collaboration 
on a related project.  We are especially grateful to Dominic Joyce for 
pointing out Lemma (\ref{fixedloclemma}), 
which resulted in significant simplifications of the 
original proofs. The work of D.-E. D. is partially supported by 
NSF grants PHY-0555374-2006 
and PHY-0854757-2009. The work of W.-y. C. was supported by DOE grant 
DE-FG02-96ER40959.

\section{Stack Function  Algebras for ADHM Quiver Sheaves}\label{stfctsect} 

This section explains how the formalism of stack functions and Ringel-Hall algebras 
constructed by Joyce in \cite{J-I}-\cite{J-IV}, \cite{J-motivic} can be applied to 
ADHM quiver sheaves on a smooth projective curve $X$ over $\IC$. 
Note that a detailed exposition of Joyce's results can be found for example in 
\cite[Sect. 2]{generating}, so we will restrict ourselves to a brief recollection of the main steps of the construction. 

\subsection{Brief Review of Joyce Theory} 
Let ${\mathfrak F}$ be an algebraic stack locally of finite type 
over $\IC$ with affine geometric stabilizers (that is, the automorphisms 
groups of $\IC$-valued points of $\mfcf$ are affine algebraic groups 
over $\IC$.) The space of stack functions of $\mfcf$ is a $\IQ$-vector 
space constructed as follows \cite[Sect. 2.3]{J-II}. 
\begin{itemize}
\item Consider pairs $({\mathfrak X}, \varrho)$ 
where ${\mathfrak X}$ is an algebraic $\IC$-stack of finite type with affine geometric 
stabilizers and $\varrho: {\mathfrak X}\to \mfcf$ is a finite type 
morphism of algebraic stacks. 
\item Two such pairs are said to be equivalent,
$({\mathfrak X}, \varrho)\sim ({\mathfrak X}', \varrho')$, if there is an 
isomorphism of stacks ${\mathfrak X}\simeq {\mathfrak X}'$
so that the obvious triangle diagram is commutative. 
 Denote equivalence classes by $[({\mathfrak X}, \varrho)]$.
\item Suppose $({\mathfrak X}, \varrho)$ is a pair as above, 
and ${\mathfrak Y} \hookrightarrow {\mathfrak X}$ is a closed substack.
Then the pair $({\mathfrak X},\varrho)$ yields two pairs 
$({\mathfrak Y}, \varrho|_{{\mathfrak Y}})$ and 
$({\mathfrak X}\setminus {{\mathfrak Y}}, \varrho|_{{\mathfrak X}\setminus {\mathfrak Y}})$. 
The stack function space ${\underline {\mathrm{SF}}}(\mfcf)$ is the $\IQ$-vector space 
generated by equivalence classes $[({\mathfrak X}, \varrho)]$ subject to the relation 
\[
[({\mathfrak X}, \varrho)] = [({\mathfrak Y}, \varrho|_{{\mathfrak Y}})] + 
[({\mathfrak X}\setminus {{\mathfrak Y}}, \varrho|_{{\mathfrak X}\setminus {\mathfrak Y}})].
\]
${\mathrm{SF}}(\mfcf)\subseteq {\underline{ \mathrm{SF}}}(\mfcf)$ is the linear subspace generated by equivalence classes of pairs 
$[({\mathfrak X}, \varrho)]$ with $\varrho$ representable. 
\end{itemize}

A central element in Joyce's theory is the existence of an associative algebra 
structure on the $\IQ$-vector space ${\mathrm{\underline {SF}}}(\mfcf)$ 
when $\mfcf$ is the moduli space of all objects in a $\IC$-linear abelian category 
$\CC$ satisfying certain assumptions \cite[Assumption 7.1]{J-I}, 
\cite[Assumption 8.1]{J-I}. 
The basic assumptions require $\CC$  to be noetherian and artinian and 
and all morphisms spaces in $\CC$  to be finite dimensional 
complex vector spaces. Natural $\IC$-bilinear composition maps 
of the form 
\[
\mathrm{Ext}^i(B,C) \times \mathrm{Ext}^j(A,B)
\to \mathrm{Ext}^{i+j}(A,C) 
\]
are required to exist for $0\leq i,j\leq 1$, $i+j=0,1$ and all $A,B,C$ objects of $\CC$. 
Moreover,  a quotient $K(\CC)$ of the Grothendieck group $K_0(\CC)$ 
by some fixed subgroup is also required, with the property 
that $[A]=0$ in $K(\CC)$ $\Rightarrow$ $A=0$ in $\CC$.
The cone spanned by classes of objects of $\CC$ in $K(\CC)$ 
will be denoted by ${\overline C}(\CC)$. The complement of the class $[0]
\in {\overline C}(\CC)$ will be denoted by $C(\CC)$. 

The remaining assumptions in \cite[Assumption 7.1]{J-I}, 
\cite[Assumption 8.1]{J-I} will not be listed in detail here. 
Essentially, one requires the existence of 
Artin moduli stacks 
$\obj(\CC)$, ${\mathfrak {Ex}}(\CC)$, locally of finite type over $\IC$, 
parameterizing all objects of $\CC$, respectively three term exact 
sequences
\be\label{eq:Aexseq} 
0\to A'\to A\to A''\to 0
\ee
in $\CC$. 
Moreover there also exist natural projections 
\be\label{eq:natstackproj} 
{\mathfrak p}, {\mathfrak p}',{\mathfrak p}'': {\mathfrak {Ex}}(\CA)
\to \obj(\CC)
\ee
which are $1$-morphisms of Artin stacks of finite type. 
There should also exist natural disjoint union decompositions 
\be\label{eq:disjunionA}
\begin{aligned}
\obj(\CC)& = \mathop{\coprod}_{\alpha \in {\overline C}(\CC)} \obj(\CC,\alpha)\\
{\mathfrak {Ex}}(\CC)& = \mathop{\coprod_{\alpha,\alpha',\alpha'' 
\in {\overline C}(\CC)}}_{\alpha=\alpha'+\alpha''} {\mathfrak {Ex}}(\CC,\alpha,\alpha',\alpha'')\\
\end{aligned}
\ee
compatible with the forgetful morphisms \eqref{eq:natstackproj}. 
All this data should satisfy additional natural conditions which will 
not be explicitly stated here. 

Granting assumptions \cite[Assumption 7.1]{J-I}, 
\cite[Assumption 8.1]{J-I}, one can define a  
$\IQ$-bilinear operation $\ast:{\underline {\mathrm{SF}}}(\obj(\CC))\times 
{\underline {\mathrm{SF}}}(\obj(\CC))\to {\underline {\mathrm{SF}}}(\obj(\CC))$ 
\cite[Def. 1]{J-II} as follows. Given two stack functions $[(\mfcx_i, {\mathfrak f}_i)]\in {\underline {\mathrm{SF}}}(\obj(\CC))$
set 
\be\label{eq:assprodA}
[(\mfcx_2, {\mathfrak f}_2)]\ast [(\mfcx_1, {\mathfrak f}_1)] = 
[(({{\mathfrak p}'}, {{\mathfrak p}''})^*({\mathfrak X}_1\times 
{\mathfrak X}_2), {\mathfrak p}\circ {\mathfrak f})].
\ee
where the stack function in the right hand side of equation 
\eqref{eq:assprodA} is determined by the following diagram
\be\label{eq:starproduct} 
\xymatrix{ 
({{\mathfrak p}'}, {{\mathfrak p}''})^*({\mathfrak X}_1\times 
{\mathfrak X}_2) \ar[d] \ar[rr]^{\mathfrak f} & & {\mathfrak {Ex}}(\CC)
\ar[d]_-{({{\mathfrak p}'}, {{\mathfrak p}''})} 
\ar[r]^-{{\mathfrak p}}& 
\obj(\CC)\\
{\mathfrak X}_1\times 
{\mathfrak X}_2\ar[rr]^-{{\mathfrak f}_1\times {\mathfrak f}_2} && 
\obj(\CC)\times \obj(\CC)& \\}
\ee
According to \cite[Thm. 5.2]{J-II}, 
$({\underline {\mathrm{SF}}}(\obj(\CC)), \ast, {\delta}_{[0]})$ 
is an associative algebra with unity, where $\delta_{[0]} = [(\mathrm{Spec}(\IC),0)]$ 
is the stack function determined by the zero object in $\CC$. 

For further reference, note that the construction of the associative stack function 
algebra can be also applied with no modification to an exact subcategory $\CA$ 
of $\CC$ (assuming that $\CC$ satisfies the above assumptions.)

Note also that 
an important element in the proof of wallcrossing formulas will be a refinement 
of the stack function algebra, 
the Ringel-Hall Lie algebra ${{\rm{SF}}}^{\sf ind}_{\sf al}(\obj(\CA))$. 
This is a Lie algebra over $\IQ$ whose underling vector space is the linear subspace 
of the stack function algebra 
spanned by stack functions with algebra stabilizers supported on virtually indecomposable 
objects. We will not review all the relevant definitions here since they will not be 
needed in the rest of the paper. We 
refer to \cite[Sect 5.1]{J-II}, \cite[Sect. 5.2]{J-II} for details. The important result for us \cite[Thm. 5.17]{J-II}
is that this linear subspace is closed under the  Lie bracket determined by the 
associative product $\ast$ , therefore it has a Lie algebra structure.

\subsection{Application to ADHM Quiver Sheaves}\label{ADHMstfct}
Let $X$ be a smooth projective curve over $\IC$. Let $M_1,M_2$ 
be fixed line bundles on $X$ equipped with a fixed isomorphism 
$M_1\otimes_X M_2{\buildrel \sim \over \longto} K_X^{-1}$. 
Recall that an abelian subcategory $\CC_X$ of ADHM quiver sheaves with 
twisting data $(M_1,M_2)$ 
has been defined in \cite[Sect. 3.1]{chamberI}. 
For completness, recall that the objects of $\CC_X$ are ADHM 
quiver sheaves on $X$ with $E_\infty=V\otimes \CO_X$, where 
$V$ is a finite dimensional complex vector space. Morphisms are 
natural morphisms of ADHM quiver sheaves with component at $\infty$ 
of the form $f\otimes 1_{\CO_X}$, where $f$ is a $\IC$-linear map. 

Since the objects of $\CC_X$ are decorated pairs of coherent $\CO_X$-modules, 
the basic assumptions recalled in the previous section hold for $\CC_X$. 
The quotient $K(\CC_X)$ of the Grothendieck group of $\CC_X$ is 
isomorphic to the lattice $\IZ^{3}$. The class of an 
object $\CE$ of $\CC_X$ is given by the triple $(r,e,v)=(r(\CE), d(\CE), v(\CE))
\in \IZ_{\geq 0}\times 
\IZ \times \IZ_{\geq 0}$, where $r(\CE)$, $d(\CE)$ are the rank, respectively degree 
of the underlying $\CO_X$-module $E$, and $v(\CE)$ is the dimension of $V$.

Let $\CA_X$ be the exact 
full subcategory of $\CC_X$  consisting of locally free ADHM  quiver sheaves on $X$. 
According to \cite[Lemma 5.2]{chamberI}, there is a locally finite type 
algebraic moduli stack $\obj(\CX)$ with affine geometric stabilizers 
parameterizing all objects of $\CA_X$. 
Moreover, \cite[Lemma 5.2]{chamberI}, 
there also exists an algebraic moduli stack ${\mathfrak {Ex}}(\CX)$ 
of three term exact sequences of objects of $\CA_X$, 
which is locally of finite type over $\IC$. 

Let  $\mfcf = \obj(\CX)$ 
in the construction described in the previous section.
The  remaining conditions in \cite[Assumption 7.1]{J-I}
follow by analogy with \cite[Thm. 10.10]{J-I}, \cite[Thm 10.12]{J-II} since 
the objects of $\CA_X$ are decorated sheaves on $X$. 
In conclusion,  the construction of the associative product in
\cite[Def. 5.1]{J-II} carries over to the 
present situation. Therefore we obtain again an associative algebra with unity 
 $({\underline{\mathrm{SF}}}(\obj(\CX)),\ast, \delta_{[0]})$ over $\IQ$.
 The construction of the Ringel-Hall Lie algebra of virtually indecomposable 
 representable stack functions with algebra stabilizers also carries over to the present
 case, 
 resulting in a Lie algebra $\mathrm{SF}_{\sf al}^{\sf ind}(\obj(\CX))$.

 \subsection{Stack function identities} 
 According to  \cite[Cor. 5.6]{J-I}, for any stability parameter 
$\delta \in \IR$ and any splitting type ${\sf t}\in {\sf T}$ 
there are  open immersions 
 \be\label{eq:opimmC}
\begin{aligned} 
& \obj_\delta^{ss}(\CX,r,e,1)
\hookrightarrow \obj(\CX)_{\leq 1}\hookrightarrow \obj(\CX) \\
&
\obj_\delta^{ss}(\CX,r,e,0) \hookrightarrow 
\obj(\CX)_{\leq 1}\hookrightarrow \obj(\CX). \\
\end{aligned}
\ee
The corresponding elements of the stack function algebra 
will be denoted by 
\[
\mfd_\delta(\alpha), \mfh(\alpha)\in {\underline{\mathrm{SF}}}(\obj(\CX)).
\]
where $\alpha = (r,e) \in \IZ_{\geq 1}\times \IZ$. 

Let $\delta_c \in \IR_{>0}$ be a critical stability parameter for ADHM 
sheaves on $X$ of type 
$\alpha=(r,e)\in \IZ_{\geq 1}\times \IZ$. According to \cite[Lemm. 4.13]{chamberI},
any $\delta_c$-semistable 
object of $\CA_X$ with $v=1$ has a one-step Harder-Narasimhan 
filtration with respect to $\delta_\pm$ stability, 
where $\delta_- < \delta_c$, $\delta_+> \delta_c$ 
are noncritical stability parameters sufficiently close to $\delta_c$. 
More precisely, 
let $\epsilon_\pm\in \IR_{>0}$ be positive real numbers 
as in \cite[Lemm. 4.13]{chamberI},
 for $\delta_i=\delta_c$. 
Let 
$\delta_+\in (\delta_c,\ \delta_c +\epsilon_+)$, $\delta_-\in (\delta_c-\epsilon_-, 
\ \delta_c)$ be noncritical stability parameters of type $(r,e)$. 
For simplicity, the stack functions 
${\mathfrak d}_{\delta_\pm}(\alpha)$, ${\mathfrak d}_{\delta_c}(\alpha)$
will be denoted by ${\mathfrak d}_{\pm}(\alpha)$, ${\mathfrak d}_{c}(\alpha)$
respectively. Given any numerical type $\alpha=(r,e)\in\IZ_{\geq 1}\times \IZ $, set 
and 
\[
\mu(\alpha) = {e\over r}, \qquad 
\mu_\pm (\alpha) = \mu(\alpha) + {\delta_\pm \over r},\qquad 
 \mu_c(\alpha) = \mu(\alpha) + {\delta_c \over r} 
\]
provided that  $r({\sf t})\neq 0$.

Then  the following lemma holds.
\begin{lemm}\label{SFlemmA} 
The following relations hold in ${\underline{\rm{SF}}}
(\obj(\CX)_{\leq 1})$ 
\be\label{eq:SFrelA}
\begin{aligned}
\mfd_{c}(\alpha)= \sum_{\substack{\alpha_1,\alpha_2\in \IZ_{\geq 1}\times \IZ\\ 
\alpha_1+\alpha_2=\alpha \\ 
\mu_{c}(\alpha_1)=\mu(\alpha_2)=\mu_c(\alpha)\\}}
{\mathfrak h}(\alpha_2)\ast \mfd_+(\alpha_1) \qquad 
\mfd_{c}(\alpha)  = \sum_{\substack{\alpha_1,\alpha_2\in \IZ_{\geq 1}\times \IZ\\ 
\alpha_1+\alpha_2=\alpha \\ 
\mu_{c}(\alpha_1)=\mu(\alpha_2)=\mu_c(\alpha)\\}}
\mfd_-(\alpha_1)\ast{\mathfrak h}(\alpha_2)
\end{aligned}
\ee
where the sums in the right hand sides of equations \eqref{eq:SFrelA} are finite. 
\end{lemm}

{\it Proof.} 
Given \cite[Lemm. 4.13]{chamberI}, theorem \cite[Thm 5.11]{J-IV} 
applies to the present case, yielding formulas \eqref{eq:SFrelA}.
 Finiteness is obvious from the summation conditions. 

\hfill$\Box$

The important point in the following is that relations \eqref{eq:SFrelA} can be 
inverted according to \cite[Thm 5.12]{J-IV}. 
In order to write down the inverse relations, for any $l\in \IZ_{\geq 1}$ and 
any $1\leq j\leq l$ define 
\be\label{eq:ljsplit}
\begin{aligned}
{\sf S}_{\delta_c}^{(l,j)}(\alpha) = 
\bigg\{ &(\alpha_1,\ldots,\alpha_l)\in (\IZ_{\geq 1}\times\IZ)^{\times l}\, \bigg|\, 
 \sum_{i=1}^l\alpha_i =\alpha,\\
& \mu(\alpha_i)= \mu_c(\alpha),\ 1\leq i\leq l, \
i\neq j,\ \mu_c(\alpha_j)=\mu_c(\alpha)\bigg\}\\
\end{aligned}
\ee
Obviously, ${\sf S}_{\delta_c}^{(l,j)}(\alpha)$ is a finite set 
for fixed $\alpha, l, j$. 

Then \cite[Thm 5.12]{J-IV} 
implies the following 
\begin{lemm}\label{SFlemmB} 
The following relations hold in ${\underline{\rm{SF}}}
(\obj(\CX))$ 
\be\label{eq:SFrelB} 
\begin{aligned} 
\mfd_{+}(\alpha) & = \mfd_c(\alpha)+ \mathop{\sum_{l\geq 2}}_{} (-1)^{l-1} 
\mathop{\sum_{
({\alpha}_1,\ldots, {\alpha}_l)\in {\sf S}_{\delta_c}^{(l,l)}(\alpha)}}_{} 
 {\mathfrak h}({\alpha}_1) \ast\cdots \ast 
{\mathfrak h}({\alpha}_{l-1})\ast \mfd_c({\alpha}_l)\\
\mfd_{-}({\alpha}) & = \mfd_c({\alpha})+
\mathop{\sum_{l\geq 2}}_{} (-1)^{l-1} \mathop{\sum_{
({\alpha}_1,\ldots, {\alpha}_l)\in {\sf S}_{\delta_c}^{(l,1)}(\alpha)}}_{} 
\mfd_c({\alpha}_1)\ast {\mathfrak h}({\alpha}_2)\ast\cdots \ast 
{\mathfrak h}({\alpha}_{l})\\
\end{aligned} 
\ee
where the sums in the right hand sides of equations \eqref{eq:SFrelB} 
are finite. 
\end{lemm}

{\it Proof.} We will check only the first equation in 
\eqref{eq:SFrelB} since the second is entirely analogous. 
According to \cite[Thm 5.12]{J-IV}, inverting the first relation 
in \eqref{eq:SFrelA} yields 
\[
\begin{aligned} 
\mfd_{+}({\alpha}) & = \mathop{\sum_{l\geq 1}}_{} (-1)^{l-1} \mathop{\sum_{j=1}^l}_{}
\mathop{\sum_{
({\alpha}_1,\ldots, {\alpha}_l)\in {\sf S}_{\delta_c}^{(l,j)}(\alpha)}}_{\mu_+^j({\alpha}_1+\cdots
+{\alpha}_k) < \mu_+^j({\alpha}_{k+1}+\cdots + {\alpha}_l)} 
{\mathfrak h}({\alpha}_1)\ast\cdots \ast \mfd_c({\alpha}_j)\ast \cdots \ast 
{\mathfrak h}({\alpha}_l)\\
\end{aligned}
\]
where 
\[
\mu_+^{j}({\alpha}_1+\cdots+{\alpha}_k) = \left\{\begin{array}{ll} 
\mu({\alpha}_1+\cdots+{\alpha}_k) & \mathrm{for}\ k<j \\
\mu_+({\alpha}_1+\cdots+{\alpha}_k) & \mathrm{for}\ k\geq j\\
\end{array}\right.
\]
\[
\mu_+^{j}({\alpha}_{k+1}+\cdots+{\alpha}_l) = \left\{\begin{array}{ll} 
\mu({\alpha}_{k+1}\cdots+{\alpha}_l) & \mathrm{for}\ k\geq j  \\
\mu_+({\alpha}_{k+1}+\cdots+{\alpha}_l) & \mathrm{for}\ k<j\\
\end{array}\right.
\]
for any $l\geq 2$ and any $1\leq k\leq l-1$. 
However, using the 
relations 
\[ 
\mu({\alpha}_i) = \mu_c({\alpha}_j) = \mu_c({\alpha}) 
\]
in \eqref{eq:ljsplit}  and $\delta_+>\delta_c$, 
it is straightforward to prove that the inequality 
\[ 
\mu_+^j({\alpha}_1+\cdots
+{\alpha}_k) < \mu_+^j({\alpha}_{k+1}+\cdots + {\alpha}_l)
\]
is satisfied if and only if $j=l$. 

\hfill $\Box$

Lemmas (\ref{SFlemmA}), (\ref{SFlemmB}) imply the following corollary, which follows 
by direct substitution.
\begin{coro}\label{SFcor}
Under the conditions of lemmas (\ref{SFlemmA}), (\ref{SFlemmB}) the following relations 
hold in the stack function algebra $\underline{SF}(\obj(\CX))$.
\be\label{eq:SFrelC} 
\begin{aligned}
\mfd_+({\alpha})-\mfd_-({\alpha}) = \sum_{l\geq 2} (-1)^{l} 
\mathop{\sum_{
({\alpha}_1,\ldots, {\alpha}_l)\in {\sf S}_{\delta_c}^{(l,l)}(\alpha)}}_{} 
 \mfh({\alpha}_1)\ast \cdots 
\ast [\mfd_-({\alpha}_l), \mfh({\alpha}_{l-1})]
\end{aligned}
\ee
\end{coro}

Next note that since the sum in the right hand side of equation 
\eqref{eq:SFrelC} is finite, the parameter $\delta_- \in \IR_{>0}$ 
can be chosen to be noncritical with respect to all types $\alpha_l$
 so that $\mfd_-({{\alpha}_l})\neq 0$. 
This implies that any {\bf S}-fixed $\delta_-$-semistable object 
of splitting type ${\alpha}_l$ is $\delta_-$-stable. 
Since $\delta_\pm$ have been chosen noncritical of type $(r,e)$
the same holds for $\delta_\pm$-semistable objects of splitting 
type ${\alpha}$. In particular the automorphism group of all 
such objects is isomorphic to $\IC^\times$, according to 
\cite[Lemm. 3.7]{chamberI}. Given the definition of 
virtually indecomposable objects with algebra stabilizers 
\cite[Sect. 5.1-5.2]{J-II}, this implies that 
the stack functions $\mfd_\pm({\alpha})$, $\mfd_-({\alpha}_l)$ belong to the 
Lie algebra ${\mathrm{SF}}^{\sf ind}_{\sf al}(\obj(\CX))$ for all possible 
splitting types  ${\alpha}_l$ 
in the right hand side of equation 
\eqref{eq:SFrelC}. 

However, the stack functions 
${\mathfrak h}({\alpha}_i)$ in the same equation  
do not satisfy this property for arbitrary 
splitting type ${\alpha}_i$, since 
strictly semistable Higgs sheaves will be present.  Then one has to 
use \cite[Thm. 8.7]{J-III} in order to construct 
virtually indecomposable log stack functions $\mfg(\alpha)$ as follows 
\be\label{eq:logstfct}
\begin{aligned} 
{\mathfrak g}(\alpha)& =\sum_{l\geq 1} {(-1)^{l-1}\over l}
\sum_{\substack{{\alpha}_1,\ldots ,
{\alpha}_l\in \IZ_{\geq 1}\times \IZ\\
\alpha_1+\cdots +\alpha_l =\alpha\\
 \mu(\alpha_i)=\mu(\alpha), \ 1\leq i\leq l\\}}
{\mathfrak h}({\alpha}_1)\ast \cdots  \ast 
{\mathfrak h}({\alpha}_l)\\
\end{aligned}
\ee
where the sum in the right hand side is finite. 
Then \cite[Thm. 8.7]{J-III} implies that $\mfg({\alpha})$ is an element of the 
Lie algebra ${\mathrm{SF}}^{\sf ind}_{\sf al}(\obj(\CX))$. 
Moreover, the following inverse relation holds \cite[Thm 8.2]{J-III}
\be\label{eq:expstfct} 
\begin{aligned} 
{\mathfrak h}({\alpha} )& =\sum_{l\geq 1} {1\over l!}
\sum_{\substack{{\alpha}_1,\ldots ,
{\alpha}_l\in \IZ_{\geq 1}\times \IZ\\
\alpha_1+\cdots +\alpha_l =\alpha\\
 \mu(\alpha_i)=\mu(\alpha), \ 1\leq i\leq l\\}}
{\mathfrak g}({\alpha}_1)\ast \cdots  \ast 
{\mathfrak g}({\alpha}_l)\\
\end{aligned}
\ee
where the sum in the right hand side is again finite. 
\begin{lemm}\label{SFlemmD} 
The following relation holds in ${\mathrm{SF}}^{\sf ind}_{\sf al}(\obj(\CX))$
\be\label{eq:SFrelD}
\begin{aligned}
\mfd_+({\alpha}) - \mfd_-({\alpha}) = &  \sum_{l\geq 2}{(-1)^{l-1}\over (l-1)!} 
\mathop{\sum_{
({\alpha}_1,\ldots, {\alpha}_l)\in {\sf S}_{\delta_c}^{(l,l)}(\alpha)}}_{} 
[\mfg({\alpha}_1), \cdots [\mfg({\alpha}_{l-1}),\mfd_{-}({\alpha}_l)] \cdots ]
\end{aligned}
\ee
\end{lemm}

{\it Proof.}
 Expanding the commutators in each term in the right hand side of 
equation \eqref{eq:SFrelD} yields 
\[
\begin{aligned}
&  [\mfg({\alpha}_1), \cdots [\mfg({\alpha}_{l-1}),\mfd_{-}({\alpha}_l)] \cdots ]= \\
&\sum_{k=0}^{l-1} \mathop{\sum_{i_1,\ldots, i_k=1}^{l-1}}_{i_1<\cdots <i_k} 
\ \ \mathop{\sum_{j_1,\ldots,j_{l-1-k}\in \{2,\ldots, l\}\setminus
\{i_1,\ldots, i_k\}}^l}_{j_1<\cdots<j_{l-1-k}} (-1)^{k} 
\mfg({\alpha}_{i_1})\ast \cdots \ast\mfg({\alpha}_{i_k}) 
\ast \mfd_-({\alpha}_l)\ast 
\mfg({\alpha}_{j_{l-1-k}})\ast \cdots \mfg({\alpha}_{j_1})\\
\end{aligned}
\]
where, by convention,  $\{i_1,\ldots, i_k\}=\emptyset$, $\{j_1,\ldots, j_{l-1-k}\} = 
\{1,\ldots, l-1\}$ if $k=0$, respectively $\{i_1,\ldots, i_k\}=\{1,\ldots, l-1\}$, $\{j_1,\ldots, j_{l-1-k}\} = \emptyset$ if $k=l-1$. 
Summing over all values of 
$({\alpha}_1, \ldots, {\alpha}_l)\in {\sf S}_{\delta_c}^{(l,l)}(\alpha)$  for fixed $l\geq 2$
yields
\be\label{eq:rhA}
{(-1)^{l-1}\over (l-1)!}\mathop{\sum_{
({\alpha}_1,\ldots, {\alpha}_l)\in {\sf S}_{\delta_c}^{(l,l)}(\alpha)}}_{} 
\sum_{k=0}^{l-1} (-1)^{k} {l-1 \choose k} 
\mfg({\alpha}_{1})\ast \cdots \ast\mfg({\alpha}_{k}) 
\ast \mfd_-({\alpha}_l) 
\ast \mfg({\alpha}_{k+1})\ast \cdots \ast 
\mfg({\alpha}_{{l-1}}).
\ee
employing similar conventions. 
Substituting \eqref{eq:expstfct} in \eqref{eq:SFrelC}, we obtain 
\be\label{eq:lhA}
\begin{aligned} 
& \mfd_+({\alpha})-\mfd_-({\alpha}) =  \\
& \sum_{p\geq 2} (-1)^{p} \mathop{\sum_{
({\alpha}_1,\ldots, {\alpha}_p)\in {\sf S}_{\delta_c}^{(p,p)}(\alpha)}}_{} 
\mathop{\sum_{m_{1}\geq 1}}_{} \mathop{\sum_{({\beta}_{1,1}, \ldots, 
{\beta}_{1,m_1})\in {\sf S}^{(m_1)}({\alpha}_1)}}_{} \cdots 
\mathop{\sum_{m_{p-1}\geq 1}}_{} \mathop{\sum_{({\beta}_{p-1,1}, \ldots ,
{\beta}_{p-1,m_{p-1}})\in {\sf S}^{(m_{p-1})}({\alpha}_{p-1})}}_{}\\
& {1\over m_1! \cdots m_{p-1}!} 
\mfg({\beta}_{1,1})\ast \cdots 
\ast \mfg({\beta}_{1,m_1})\ast \mfg({\beta}_{2,1})\ast \cdots 
\ast \mfg({\beta}_{2,m_2})\ast \cdots \\
& \qquad \qquad \qquad \qquad \qquad \qquad \qquad \qquad \qquad \ast 
[ \mfd_-({\beta}_p),\, \mfg({\beta}_{p-1,1})\ast \cdots 
\ast \mfg({\beta}_{p-1,m_{p-1}})]\\
 \end{aligned}
\ee
where 
\[
{\sf S}^{(m)}({\alpha})=\big\{ (\beta_1,\ldots,\beta_m) \in (\IZ_{\geq 1}\times \IZ)^m\, \big|\, 
\beta_1+\cdots + \beta_m =\alpha,\ \mu(\beta_1)=\cdots = \mu(\beta_m) \big\}
\]
for any $m\geq 1$ and any $\alpha \in \IZ_{\geq 1}\times \IZ$. 

The right hand side of \eqref{eq:lhA} can be rewritten as 
\be\label{eq:lhB} 
\begin{aligned} 
& \mfd_+({\alpha})-\mfd_-({\alpha}) =  \\
& \sum_{p\geq 2} (-1)^{p} 
\mathop{\sum_{m_1,\ldots,m_{p-1} \geq 1}}_{}
\mathop{\sum_{({\beta}_{1,1},\ldots, {\beta}_{1,m_1},\ldots, 
{\beta}_{p-1,1},\ldots, {\beta}_{p-1,m_{p-1}},{\beta}_l)\in 
{\sf S}^{(l,l)}_{\delta_c}(\alpha)}}_{}
 {1\over m_1! \cdots m_{p-1}!} \\
 &
\mfg({\beta}_{1,1})\ast \cdots 
\ast \mfg({\beta}_{1,m_1})\ast \mfg({\beta}_{2,1})\ast \cdots 
\ast \mfg({\beta}_{2,m_2})\ast \cdots \ast 
[ \mfd_-({\beta}_l),\, \mfg({\beta}_{p-1,1})\ast \cdots 
\ast \mfg({\beta}_{p-1,m_{p-1}})]\\
 \end{aligned}
\ee
where $l =m_1+\cdots + m_{p-1} +1$. 

Note that for fixed $(p,l)$ 
in the right hand side of \eqref{eq:lhB} we sum over 
ordered sequences $(m_1,\ldots,m_{p-1})\in \IZ_{>0}^{p-1}$ satisfying $m_1+\cdots + m_{p-1}=l-1$. For $p\geq 3$ 
there are exactly two monomials associated to each such ordered sequence, 
namely
\[ 
 \mfg({\beta}_{1})\ast \cdots \ast 
\mfg({\beta}_{l-1})\ast \mfd_-({\beta}_l) \qquad \mathrm{and}\qquad 
\mfg({\beta}_{1})\ast \cdots \ast\mfg({\beta}_{k}) 
\ast \mfd_-({\beta}_l) 
\ast \mfg({\beta}_{k+1})\ast \cdots \ast 
\mfg({\beta}_{l-1})
\] 
with $1\leq k =m_{p-2} \leq l-1$. The same statement holds for $p=2$, 
except that the second monomial in the above equation reads 
$\mfd_-({\beta}_l) 
\ast \mfg({\beta}_{k+1})\ast \cdots \ast 
\mfg({\beta}_{l-1})$.

Given an arbitrary monomial of the form 
\be\label{eq:monomA} 
 \mfg({\beta}_{1})\ast \cdots \ast 
\mfg({\beta}_{l-1})\ast \mfd_-({\beta}_l)
 \ee
with fixed $l\geq 2$ and fixed $({\sf t}_1,\ldots, {\sf t}_l)\in 
{\sf S}^{(l,l)}_{\delta_c}({\alpha})$ 
there is an obvious one-to-one correspondence 
between ordered sequences $(m_1, \ldots, m_{p-1})$ and partitions 
of the ordered sequence $({\beta}_1,\ldots, {\beta}_{l-1})$ \
of the form 
\be\label{eq:monompartitionA}
\left({{\beta}_{1}}, \ldots, {{\beta}_{m_1}}\ |\
\ldots \ |\ {{\beta}_{l-m_p}}\ldots {{\beta}_{l-1}}\ \right).
\ee
Moreover, the sequence $(m_1,\ldots, m_{p-1})$ also determines a length $(p-1)$ 
unordered partition 
$\lambda_{(m_1,\ldots, m_{p-1})} =(1^{j_1}, \ldots, s^{j_s})$
of $(l-1)$, which will be called the underlying partition 
of the sequence $(m_1,\ldots, m_{p-1})$. The factor  
\[
{1\over m_1! \cdots m_{p-1}!} = {1\over (1!)^{j_1} \cdots (s!)^{j_s}}
\]
depends only on the underlying partition $\lambda_{(m_1,\ldots, m_{p-1})}$. 
Conversely, for  a fixed length $(p-1)$ partition 
$\lambda=(1^{j_1}, 2^{j_2}, \ldots, s^{j_s})$ of $(l-1)$, 
there are
\[
{(p-1)!\over j_1! j_2! \cdots  j_s!}
\]
distinct ordered sequences $(m_1,\ldots, m_{p-1})$ as above with 
underlying partition $\lambda$. Each such sequence corresponds to a 
partition of the set $({\beta}_1,\ldots, {\beta}_{l-1})$ of the form \eqref{eq:monompartitionA}. 

Similar arguments apply to any monomial of the form 
\be\label{eq:monomB} 
\mfg({\beta}_{1})\ast \cdots \ast\mfg({\beta}_{k}) 
\ast \mfd_-({\beta}_l) 
\ast \mfg({\beta}_{k+1})\ast \cdots \ast 
\mfg({\beta}_{l-1})
\ee
with $1\leq k\leq l-1$. For $p\geq 3$, there is a one-to-one correspondence between 
ordered sequences $(m_1,\ldots, m_{p-2})$ with $m_{p-2}=k$ and partitions of 
the ordered sequence $({\beta}_{1}, \ldots, {\beta}_k)$ of the form 
\be\label{eq:monompartitionB}
\left({{\beta}_{1}}, \ldots, {{\beta}_{m_1}}\ |\
\ldots \ |\ {{\beta}_{m_{p-3}+1}}\ldots {{\beta}_{k}}\ \right)
\ee
Moreover, an ordered sequence $(m_1,\ldots, m_{p-2})$ as above also determines 
a length $(p-2)$ partition of $k$, $\lambda_{(m_1,\ldots, m_{p-2})}=(1^{s_1},\ldots, 
s^{j_s})$.  The following relation holds 
\[
{1\over m_1!\cdots m_{p-2}!} = {1\over (1!)^{j_1} \cdots (s!)^{j_s}}.
\]
Conversely, for 
a length $(p-2)$ partition of $k$, $\lambda=(1^{s_1},\ldots, 
s^{j_s})$ there are
\[
{ (p-2)! \over j_1! \cdots j_s!}
\]
distinct ordered sequences $(m_1,\ldots, m_{p-2})$ with underlying partition 
$\lambda$. 

In conclusion, 
the right hand side of \eqref{eq:lhB} can be further rewritten as follows 
\be\label{eq:lhC} 
\begin{aligned} 
& \mfd_+({\alpha})-\mfd_-({\alpha}) =  \\ & \sum_{l\geq 2} 
\mathop{\sum_{({\beta}_1,{\beta}_2,\ldots, {\beta}_l)\in {\beta}^{(l,l)}_{\delta_c,{\beta}}}}_{} \mathop{\sum_{k=0}^{l-1}}_{} 
c_k({\beta}_1, \ldots, {\beta}_l)  \mfg({\beta}_{1})\ast \cdots \ast\mfg({\beta}_{k}) 
\ast \mfd_-({\beta}_l) 
\ast \mfg({\beta}_{k+1})\ast \cdots \ast 
\mfg({\beta}_{l-1})\\
\end{aligned}
\ee
 where the coefficients $c_k({\beta}_1, \ldots, {\beta}_l)$ are given by 
 \[
 \begin{aligned} 
 c_{l-1}({\beta}_1, \ldots, {\beta}_l) & = -\mathop{\sum_{p\geq 2}}_{} 
 (-1)^{p} \mathop{\sum_{\lambda \in {\sf P}_{p-1}(l-1)}}_{
\lambda = (1^{j_1}, 2^{j_2}, \ldots, s^{j_s})}{(p-1)!\over j_1! j_2! \cdots  j_s!}
{1\over (1!)^{j_1} \cdots (s!)^{j_s}}\\
\end{aligned}
\]
\[
\begin{aligned} 
 c_k({\beta}_1, \ldots, {\beta}_l) & = {1\over (l-k-1)!} 
 \mathop{\sum_{p\geq 3}}_{} 
 (-1)^{p} \mathop{\sum_{\lambda \in {\sf P}_{p-2}(k)}}_{
\lambda = (1^{j_1}, 2^{j_2}, \ldots, s^{j_s})}{(p-2)!\over j_1! j_2! \cdots  j_s!}
{1\over (1!)^{j_1} \cdots (s!)^{j_s}}\\
\end{aligned}
\]
for $1\leq k\leq l-2$, $l\geq 3$, and 
\[
c_0({\beta}_1, \ldots, {\beta}_l) = {1\over (l-1)!} 
\]
if $k=0$. 
 Here, 
where ${\sf P}_{p-1}(l-1)$ denotes the set of length $(p-1)$ 
partitions
of $(l-1)$, ${\sf P}_{p-2}(k)$ denotes the set of length $(p-2)$ 
partitions
of $k$.

Next note that the coefficients $c_k({\beta}_1, \ldots, {\beta}_l)$
may be expressed in terms of Bell polynomials 
\[
\begin{aligned}
c_{l-1}({\beta}_1, \ldots, {\beta}_l) & = {1\over (l-1)!} \mathop{\sum_{p\geq 2}}_{} 
 (-1)^{p-1} (p-1)! B_{l-1,p-1}(1,1,\ldots, 1) \\
 \end{aligned}
 \]
 respectively
 \[
 \begin{aligned}
   c_k({\beta}_1, \ldots, {\beta}_l) & = {1\over k!\, (l-k-1)!} 
 \mathop{\sum_{p\geq 2}}_{}(-1)^{p-2} (p-2)! B_{l-k-1,p-2}(1,1,\ldots,1).\\
 \end{aligned}
 \] 
 for $1\leq k \leq l-2$, $l\geq 3$. 
 Some basic facts on Bell polynomials are recalled fro convenience in 
 appendix A. 
  Then a special case of the Fa{\`a} di Bruno formula 
 (see equation \eqref{eq:Bellsum}) yields 
 \[
 c_{l-1}({\beta}_1, \ldots, {\beta}_l) = {(-1)^{l-1}\over (l-1)!} \qquad 
  c_k({\beta}_1, \ldots, {\beta}_l) = {(-1)^{l-k-1}\over k!\, (l-k-1)!}= {(-1)^{l-k}\over (l-1)!} 
  {l-k-1 \choose k}.
  \]
  Therefore, taking into account equation \eqref{eq:rhA}, 
  the final formula for the difference $\mfd_+({\alpha}) - \mfd_-(\alpha)$ 
  is indeed \eqref{eq:SFrelD}. 
 
 \hfill $\Box$ 
 
 Analogous arguments yield an identity relating the stack functions 
$\mfd_\pm(\alpha)$ where $\delta_+\in \IR_{>0}$, $\delta_-\in \IR_{<0}$ are stability 
parameters sufficiently close to the origin. 
More precisely, take $\delta_+<\epsilon_+$, $\delta_->\epsilon_-$, 
where $\epsilon_\pm$ are as in \cite[Lemm. 4.15]{chamberI}.
Let ${\mathfrak o}$ be the stack function determined by the moduli stack of objects 
of $\CA_X$ of type $(r,e,v)=(0,0,1)$. Note that any such object is isomorphic to 
$O=(0,\IC,0,0,0,0)$ and the moduli stack in question is isomorphic to 
the quotient stack $[*/\IC^\times]$. Let ${\sf S}^{(l)}_{0}(\alpha)$ be the set 
obtained by setting $\delta_c=0$ in equation  
in \eqref{eq:ljsplit}, which then becomes independent of $1\leq j\leq l$.  
Then, in complete analogy with lemmas (\ref{SFlemmA}), (\ref{SFlemmB}), 
(\ref{SFlemmD}),
\begin{lemm}\label{SFlemmE} 
The following identity holds in the Ringel-Hall Lie algebra 
${\mathrm{SF}}_{\sf al}^{\sf ind}(\obj(\CX))$
\be\label{eq:SFrelF} 
\begin{aligned}
\mfd_+({\alpha}) -\mfd_-({\alpha}) = & \mathop{\sum_{l\geq 2}}_{} 
{(-1)^{l-1}\over (l-1)!} \mathop{\sum_{({\alpha}_1,\ldots, {\alpha}_l)\in 
{\sf S}^{(l)}_{0}({\alpha})}}_{}[\mfg({\alpha}_1),\cdots 
[\mfg({\alpha}_{l-1}),\mfd_-({\alpha})_l]\cdots ]\\
& +\mathop{\sum_{l\geq 1}}_{} 
{(-1)^{l}\over l!} \mathop{\sum_{({\alpha}_1,\ldots, {\alpha}_l)\in 
{\sf S}^{(l)}_{0}({\alpha})}}_{}[\mfg({\alpha}_1),\cdots [\mfg({\alpha}_l),{\mathfrak o}]\cdots ]\\
\end{aligned}
\ee
where the sum in the right hand side of equation \eqref{eq:SFrelF} is finite. 
\end{lemm}
  
\section{Wallcrossing Formulas}\label{crtwallsect} 

In this section we prove  theorems (\ref{wallcrossingthmB}) and 
(\ref{wallcrossingthmC}).  

\subsection{ADHM invariants via Behrend's constructible function}
In order to derive wallcrossing formulas using the formalism of Joyce and Song, 
the ADHM invariants must be first expressed in terms of Behrend's weighted Euler characteristic. This is the content of the following lemma, which was pointed out by  Dominic 
Joyce.

\begin{lemm}\label{fixedloclemma} 
Let $\delta\in \IR_{>0}$ be a noncritical stability parameter of type 
$(r,e)$. Then 
\be\label{eq:adhminvA}
A_\delta^{\bf S}(r,e) =  \chi^{B}(\mfm_\delta^{ss}(\CX,r,e)) 
\ee
where the right hand side of equation \eqref{eq:adhminvA} is Behrend's weighted Euler 
characteristic of the algebraic space $\mfm_\delta^{ss}(\CX,r,e)$.
 \end{lemm} 

{\it Proof.} 
Recall that the ADHM invariant $A_\delta^{\bf S}(r,e)$ is defined by virtual 
integration on the fixed locus $\mfm_\delta(\CX,r,e)^{\bf S}$
\[
A_\delta^{\bf S}(r,e) = \int_{[\mfm_\delta^{ss}(\CX,r,e)^{\bf S}]} 
e_{\bf S}(N^{vir}_{\mfm_\delta^{ss}(\CX,r,e)^{\bf S}/\mfm_\delta^{ss}(\CX,r,e)})^{-1} 
\]
The virtual cycle of the fixed locus is determined by the fixed part of the perfect tangent-obstruction theory of the moduli space restricted to the fixed locus. 
Since the perfect tangent-obstruction theory of $\mfm_\delta(\CX,r,e)$ is {\bf S}-equivariant 
symmetric it follows that the induced tangent-obstruction theory of the fixed locus is 
 symmetric. Therefore the resulting virtual cycle 
is a $0$-cycle. 

The virtual normal bundle $N^{vir}_{\mfm_\delta(\CX,r,e)^{{\bf S}}/\mfm_
\delta(\CX,r,e)}$ 
is determined by the {\bf S}-moving part of the perfect obstruction theory 
of $\mfm_\delta^{ss}(\CX,r,e)$ restricted to $\mfm_\delta(\CX,r,e)^{{\bf S}}$, 
which is also {\bf S}-equivariant symmetric. By construction, the virtual normal bundle 
 $N^{vir}_{\mfm_\delta(\CX,r,e)^{{\bf S}}/\mfm_\delta(\CX,r,e)}$ is an equivariant K-theory 
 class of the form 
 \be\label{eq:virnormA}
 {\IE_1^m}-{\IE_2^m},
 \ee
 where $\IE_1^m$ is an equivariant locally free sheaf on $\mfm_\delta(\CX,r,e)^{{\bf S}}$
 and $\IE_2^m$ is its (equivariant) dual. Moreover, the character decomposition of 
 $\IE_1^m$ does not contain the trivial character. 
 
 Since the virtual cycle of the fixed locus is a 0-cycle, it suffices to compute 
the equivariant Euler class $e_{\bf S}(N^{vir}_{\mfm_\delta(\CX,r,e)^{{\bf S}}/\mfm_
\delta(\CX,r,e)}\big|_{\mathfrak m})$ 
of the restriction of the virtual normal bundle to a closed point ${\mfm}$ of the fixed locus. 
Let $\CE$ be the {\bf S}-fixed $\delta$-stable ADHM sheaf on $X$ corresponding to $\mfm$. 
Then, given the construction of the perfect tangent-obstruction theory in \cite[Sect. 5.4]{chamberI} 
it follows that 
\[
N^{vir}_{\mfm_\delta(\CX,r,e)^{{\bf S}}/\mfm_
\delta(\CX,r,e)}\big|_{\mathfrak m} = \mathrm{Ext}^1(\CE,\CE)^m - 
\mathrm{Ext}^2(\CE,\CE)^m 
\]
where $\mathrm{Ext}^k(\CE,\CE)^m$, $k=1,2$ denotes the moving 
part of the ext group $\mathrm{Ext}^k(\CE,\CE)$ in the abelian category $\CA_\CX$. 
Moreover, using \cite[Prop. 3.15]{chamberI}, it is straightforward to check that 
there is an equivariant isomorphism $\mathrm{Ext}^2(\CE,\CE)^m\simeq 
(\mathrm{Ext}^1(\CE,\CE)^m)^\vee$. This implies that 
\be\label{eq:virnormeulerA} 
e_{\bf S}(N^{vir}_{\mfm_\delta(\CX,r,e)^{{\bf S}}/\mfm_
\delta(\CX,r,e)}\big|_{\mathfrak m}) = (-1)^{\mathrm{dim} \mathrm{Ext}^1(\CE,\CE)^m}.
\ee
Since the virtual normal bundle is a K-theory class of the form \eqref{eq:virnormA}, 
the right hand side of equation \eqref{eq:virnormeulerA} must be independent of $\CE$
when $\mfm$ varies within a connected component $\Xi$ of the fixed locus.  
This can be in fact confirmed by a direct computation based on the locally free complex 
given in \cite[Prop. 3.15]{chamberI}, but the details will not be needed in the following. 
Let $\sigma(\Xi)$ denote the common value of $(-1)^{\mathrm{dim} \mathrm{Ext}^1(\CE,\CE)^m}$ 
for all closed points $\mfm\in \Xi$.  
Then we obtain 
\be\label{eq:signA}
\begin{aligned} 
A_\delta^{\bf S}(r,e)_{\sf t} & = \sum_{\Xi} \sigma(\Xi) \int_{[\Xi]^{vir} } 1= 
 \sum_{\Xi} \sigma(\Xi) \chi^B(\Xi)\\
\eal
\ee
where $\chi^B(\Xi)$ denotes the weighted Euler character of the connected component 
$\Xi$ of the fixed locus. 
Next we claim that for any $\Xi$
\be\label{eq:weighteulerA}
\sigma(\Xi) \chi^B(\Xi) = \chi(\Xi,\nu|_{\Xi}) ,
\ee
where $\nu$ is Behrend's constructible function of the moduli space $\mfm_\delta^{ss}(\CX,r,e)$. 

Let  $\CE$ 
be an {\bf S}-fixed $\delta$-stable ADHM sheaf corresponding to a closed point $\mfm\in \Xi$ as above.
 Then 
\cite[Thm. 7.1]{chamberI} implies that the moduli space $\mfm_\delta^{ss}(\CX,r,e)$ is analytically locally 
isomorphic near $\mfm$ to the critical locus of a holomorphic function $\Phi: U\to \IC$, where 
$U\subset \mathrm{Ext}^1(\CE,\CE)$ is an analytic open neighborhood of the origin. Moreover, 
given the construction in \cite[Sect. 7]{chamberI}, 
$U,\Phi$ can be naturally chosen so that $U$ is preserved by the induced ${\bf S}$-action on 
$\mathrm{Ext}^1(\CE,\CE)$, and $\Phi$ is ${\bf S}$-invariant. In particular, $\Phi$ yields 
a holomorphic function $\Phi^{\bf S}$ on the fixed locus $U^{\bf S}\subset U$ so that 
$\Xi$ is analytically locally isomorphic to the critical locus of $\Phi^{\bf S}$. 
Then 
\[
\nu([\CE]) = (-1)^{\mathrm{dim}(\mathrm{Ext}^1(\CE,\CE))} (1-\chi_{top}(MF(\Phi,0)))
\]
where $MF(\Phi,0)$ is the Milnor fiber of $\Phi$ at $0\in U$, and 
$\chi_{top}$ denotes the topological Euler characteristic. 
Furthermore 
\[
\nu_{\Xi}([\CE]) = (-1)^{\mathrm{dim}(\mathrm{Ext}^1(\CE,\CE)^f)} 
(1-\chi_{top}(MF(\Phi^{\bf S},0))).
\]
where $\nu_\Xi$ is Behrend's constructible function of the fixed locus $\Xi$, and 
\[
\chi_{top}(MF(\Phi,0)) =\chi_{top}(MF(\Phi^{\bf S},0)).
\]
Therefore 
\[
\nu([\CE]) = (-1)^{\mathrm{dim}(\mathrm{Ext}^1(\CE,\CE)^m)} \nu_{\Xi}([\CE]),
\]
which implies that 
\be\label{eq:weighteulerB} 
\sigma(\Xi)\chi^B(\Xi) = \chi(\Xi, \nu|_\Xi).
\ee
Since $\Xi$ are the connected components of the {\bf S}-fixed locus, equation 
\eqref{eq:adhminvA} then follows easily. 

\hfill $\Box$

\subsection{Counting invariants and wallcrossing}\label{wallsubsection}

Let $\CE_1,\CE_2$ be two locally free ADHM sheaves on $X$ of numerical types  
$(r_1,e_1,1)$, $(r_2,e_2,0)$. Let 
\[
\begin{aligned}
\chi(\CE_1,\CE_2) = &\ \mathrm{dim}\mathrm{Ext}^{0}(\CE_1,\CE_2)-
\mathrm{dim}\mathrm{Ext}^{1}(\CE_1,\CE_2) \\
& - \mathrm{dim}\mathrm{Ext}^{0}(\CE_2,\CE_1) + 
\mathrm{dim}\mathrm{Ext}^{1}(\CE_2,\CE_1).\
\end{aligned} 
\]
 According to \cite[Lemm. 7.3]{chamberI},  
 \[
 \chi(\CE_1,\CE_2) = e_2 - r_2(g-1) 
 \]
 depends only of the numerical types of the two objects. 
 Abusing notation, we will also denote by $\chi$ the resulting bilinear form 
 on numerical invariants.

 Now let ${\sf L}(\CX)_{\leq 1}$ be the $\IQ$-vector space spanned by the formal symbols 
 $\lambda^{\alpha}$,  
 $\lambda^{(\alpha,1)}$, $\alpha \in \IZ_{\geq 1}\times \IZ$.
 Then the following antisymmetric bilinear form 
 \be\label{eq:symbbracket} 
 \begin{aligned}
 & \\
  [\lambda^{\alpha_1}, \lambda^{\alpha_2}]_{\leq 1}   & = 0\\
 [\lambda^{(\alpha_1,1)}, \lambda^{\alpha_2}]_{\leq 1}  
 & = -   [\lambda^{\alpha_2}, \lambda^{(\alpha_1,1)}]_{\leq 1}  =  
 (-1)^{e_2-r_2(g-1)}(e_2-r_2(g-1))
 \lambda^{(\alpha_1+\alpha_2,1)}\\
 [\lambda^{(\alpha_1,1)}, \lambda^{(\alpha_2, 1)}]_{\leq 1}&=0\\
   \end{aligned}
 \ee
 defines a $\IQ$-Lie algebra structure on ${\sf L}(\CX)_{\leq 1}$. 
   
Let $\mathrm{SF}_{\sf al}^{\sf ind}(\obj(\CX))_{\leq 1}$ be the truncation of 
the $\IQ$-vector space $\mathrm{SF}_{\sf al}^{\sf ind}(\obj(\CX))$
to stack functions $[(\mfcx,\varrho)]$ so that $\varrho$ factors through
the open immersion $\obj(\CX)^{\bf S}_{\leq 1} \hookrightarrow \obj(\CX)^{\bf S}$. 
Using the Lie algebra structure $[\ ,\ ]$ on $\mathrm{SF}_{\sf al}^{\sf ind}(\obj(\CX))$, we define a truncated  Lie algebra structure $[\ ,\ ]_{\leq 1}$ on 
$\mathrm{SF}_{\sf al}^{\sf ind}(\obj(\CX)^{\bf S})_{\leq 1}$ which is equal to 
$[\ ,\ ]$ if the arguments satisfy $v_1+v_2\leq 1$ and vanishes 
identically if both arguments are elements with $v=1$. 

Now let $\nu$ denote the Behrend constructible function of the algebraic 
stack $\obj(\CX)^{\bf S}_{\leq 1}$ defined in \cite[Prop. 4.4]{genDTI}. 
Then, given \cite[Thm 7.2]{chamberI}, \cite[Lemm. 7.3]{chamberI} 
and \cite[Thm. 7.4]{chamberI},  the following theorem holds
by analogy with \cite[Thm. 5.12]{genDTI}.
\begin{theo}\label{Liemorphthm}
There exists a Lie algebra morphism 
\be\label{eq:countingA} 
 \Psi: \mathrm{SF}_{\sf al}^{\sf ind}(\obj(\CX))_{\leq 1} \to {\sf L}(\CX)_{\leq 1} 
 \ee
which maps an element of $\mathrm{SF}_{\sf al}^{\sf ind}(\obj(\CX))_{\leq 1}$ 
 of numerical  type $\alpha$, respectively $(\alpha,1)$, 
 $\alpha \in \IZ_{\geq 1}\times \IZ$ to $\IQ \lambda^{\alpha}$, respectively
  $\IQ \lambda^{(\alpha,1)}$. 

Moreover, suppose $[(\mfcx,\varrho)]$ is an element of $\mathrm{SF}_{\sf al}^{\sf ind}(\obj(\CX))_{\leq 1}$ of type $(\alpha,1)$, 
where $\mfcx\to {\sf X}$ is a $\IC^\times$-gerbe over an algebraic 
space ${\sf X}$ of finite type over $\IC$, and $\varrho: \mfcx\to \obj(\CX)_{\leq 1}^{\bf 
S}$ is an open immersion. Then 
\be\label{eq:countingB} 
\Psi([(\mfcx,\varrho)]) = -\chi^B({\sf X})\lambda^{(\alpha,1)}
\ee
were $\chi^B({\sf X})$ is Behrend's weighted Euler characteristic 
of the algebraic space ${\sf X}$.  
\end{theo}

Recall that according to \cite[Cor. 5.5]{chamberI} for any noncritical 
stability parameter of type $(r,e)\in\IZ_{\geq 1}\times \IZ$ and 
the moduli stack 
$\obj^{ss}_\delta(\CX,r,e,1)$ 
is a $\IC^\times$-gerbe over the algebraic moduli space 
$\mfm_\delta^{ss}(\CX,r,e)$ 
of $\delta$-semistable ADHM sheaves of type $(r,e)$. 
Then theorem (\ref{Liemorphthm}) lemma (\ref{fixedloclemma})  imply
\begin{coro}\label{adhminvcor} 
Let $\delta\in \IR_{>0}$ be a noncritical stability parameter of type 
$(r,e)\in \IZ_{\geq 1}\times \IZ$. Then 
\be\label{eq:adhminvD}
\Psi(\mfd_\delta(\alpha))=-A^{\bf S}_\delta(r,e)\lambda^{(\alpha,1)}.
\ee
\end{coro}

In order to formulate a wallcrossing result for ADHM invariants, 
one has to also define Higgs sheaf invariants 
by 
\be\label{eq:higgsinvA}
\Psi(\mfg(\alpha))= -H(r,e) \lambda^\alpha
\ee
for any $\alpha=(r,e)\in \IZ_{\geq 1}\times\IZ$. 

By analogy with \cite{genDTI}, define the invariants ${\overline H}^{\bf S}(r,e)$ 
by the multicover formula
\be\label{eq:higgsinvB}
{ H}(r,e) = \mathop{\sum_{m\geq 1}}_{m|r,\ m|e}
{1\over m^2} {\overline H}(r/m, e/m).
\ee
Conjecturally, ${\overline H}(r,e)$ are $\IZ$-valued invariants.

 {\it Proof of Theorem (\ref{wallcrossingthmB}).}
 Formulas \eqref{eq:wallformulaAA} and \eqref{eq:wallformulaAB} 
follow by a simple computation applying the Lie algebra 
morphism $\Psi$ of theorem (\ref{Liemorphthm}) to the stack function 
identities derived in  lemmas (\ref{SFlemmD}), respectively (\ref{SFlemmE}).

\hfill $\Box$

\section{Applications} 

\subsection{Comparison with Kontsevich-Soibelman formula}\label{KSsect}

In this section we specialize the wallcrossing formula of Kontsevich and 
Soibelman \cite{wallcrossing} to ADHM invariants, and prove that it implies equation 
\eqref{eq:wallformulaAA}. 
Recall that locally free ADHM quiver sheaves on $X$ have a numerical invariants 
of the form $(r,e,v)\in \IZ_{\geq 0} \times \IZ \times \IZ_{\geq 0}$. 
The pair $(r,e)$ is denoted by ${\alpha}$ in theorem (\ref{wallcrossingthmB}). 
Let  ${\sf e}_{\alpha}=\lambda^\alpha$, ${\sf f}_{\alpha}=\lambda^{(\alpha,1)}$, 
${\alpha}\in \IZ_{\geq 1} \times \IZ$ 
be alternative notation for the generators of the Lie algebra ${\sf L}(\CX)_{\geq 1}$. 
Therefore 
\be\label{eq:chargeliebracket}
\begin{aligned}
& [{\sf e}_{{\alpha}_1}, {\sf e}_{{\alpha}_2}]_{\leq 1} = 0\\
&[{\sf f}_{{\alpha}_1}, {\sf f}_{{\alpha}_2}]_{\leq 1} = 0\\
&[{\sf f}_{{\alpha}_1}, {\sf e}_{{\alpha}_2}]_{\leq 1} =\chi(\alpha_1,\alpha_2)
 {\sf f}_{{\alpha}_1+{\alpha}_2}\\
\end{aligned}
\ee
where $\chi(\alpha_1,\alpha_2)=(-1)^{e_2-r_2(g-1)}(e_2-r_2(g-1))$. 

Let $\delta_c \in \IR_{>0}$ be a critical stability parameter of type $(r,e)\in 
\IZ_{\geq 1}\times \IZ$ as in theorem (\ref{wallcrossingthmB}). 
Then there exist $\alpha,\beta\in \IZ_{\geq 1}\times \IZ$, with 
\be\label{eq:slopecondM}
\mu_{c}(\alpha)=\mu(\beta)=\mu_c({\alpha})
\ee
so that any $\eta\in \IZ_{\geq 1}\times \IZ$ with 
\[
\mu_c(\eta) =\mu_c({\alpha})
\]
is uniquely written as 
\[
\eta = \alpha + q\beta, \qquad q\in \IZ_{\geq 0}
\]
and any $\rho\in \IZ_{\geq 1}\times \IZ$ with
\[
\mu(\rho) = \mu_c({\alpha}) 
\]
is uniquely written as 
\[
\rho = q \beta, \qquad q \in \IZ_{\geq 0}.
\]
Therefore $\alpha$ and $\beta$ generate a subcone of $\IZ_{\geq 1}\times \IZ$ 
consisting of elements of $\delta_c$-slope equal to $\mu_c({\alpha})$. 

For any $q\in \IZ_{\geq 0}$ define 
to be the following formal expressions
\begin{equation}
U_{\alpha+q\beta} = \mathrm{exp}({f_{\alpha+q\beta}})\qquad 
U_{q\beta} = \text{exp} (\sum_{m\geq1} \frac{e_{mq\beta}}{m^2})
\end{equation}
In this context, the wallcrossing formula of Kontsevich and Soibelman \cite{wallcrossing}
reads
\begin{equation}\label{eq:KSformA}
\prod_{{q\geq 0},\ q\uparrow} U_{q\beta}^{{\overline H}(q\beta)}
\prod_{{q\geq 0},\ q\uparrow} U_{\alpha + q \beta}^{A^{\bf S}_{+}(\alpha+q\beta)} = \prod_{{q\geq 0},\ q\uparrow} U_{\alpha + q \beta}^{A^{\bf S}_{-}(\alpha+q\beta)}\prod_{q\geq 0,\ q\downarrow} U_{q\beta}^{{\overline H}(q\beta)},
\end{equation}
where an up, respectively down arrow means that the factors in the corresponding 
product are taken in increasing, respectively decreasing order of $q$. 

In the following we will prove that equation \eqref{eq:KSformA} implies 
the wallcrossing formula \eqref{eq:wallformulaAA}. 
First note that given equation \eqref{eq:chargeliebracket}, 
the formal operators $U$  commute within each product over $q$ 
in equation \eqref{eq:KSformA}. Therefore \eqref{eq:KSformA}
can be rewritten as 
\[
\mathrm{exp}\big(\sum_{m\geq 1} \sum_{q\geq 0} 
{\overline H}(mq\beta){{\sf e}_{mq\beta}\over m^2}\big)
\prod_{q\geq 0}U_{\alpha + q \beta}^{A^{\bf S}_{+}(\alpha+q\beta)} 
= \prod_{{q\geq 0}} 
U_{\alpha + q \beta}^{A^{\bf S}_{-}(\alpha+q\beta)} 
\mathrm{exp}\big(\sum_{m\geq 1} \sum_{q\geq 0} 
{\overline H}(mq\beta){{\sf e}_{mq\beta}\over m^2}\big)
\]
This formula can be rewritten in terms of the rational invariants 
$H({\alpha})$ using  \eqref{eq:higgsinvB}. We obtain
\be\label{eq:KSformB}
\mathrm{exp}\big( \sum_{q\geq 0} 
{ H}(q\beta){{\sf e}_{q\beta}}\big)
\prod_{q\geq 0}U_{\alpha + q \beta}^{A^{\bf S}_{+}(\alpha+q\beta)} 
= \prod_{{q\geq 0}} 
U_{\alpha + q \beta}^{A^{\bf S}_{-}(\alpha+q\beta)} 
\mathrm{exp}\big( \sum_{q\geq 0} 
{ H}(q\beta){{\sf e}_{q\beta}}\big)
\ee
Let us denote by 
\[
\IH = \sum_{q\geq 0} 
{ H}(q\beta){{\sf e}_{q\beta}}.
\]
Therefore we obtain 
\be\label{eq:KSformC}
\prod_{q\geq 0}U_{\alpha + q \beta}^{A^{\bf S}_{+}(\alpha+q\beta)}  = 
\mathrm{exp}(-\IH) \prod_{{q\geq 0}} 
U_{\alpha + q \beta}^{A^{\bf S}_{-}(\alpha+q\beta)} \mathrm{exp}(\IH).
\ee
Using again the Lie algebra structure \eqref{eq:chargeliebracket}, 
note that 
\[
\prod_{q\geq 0} U_{\alpha + q \beta}^{A^{\bf S}_{\pm}(\alpha+q\beta)} = \mathrm{exp} 
\big(\sum_{q\geq 0} A^{\bf S}_{\pm}(\alpha+q\beta) {\sf f}_{\alpha+q\beta} \big)
\]
Therefore equation \eqref{eq:KSformC}) simplifies to 
\be\label{eq:KSformD}
\mathrm{exp} 
\big(\sum_{q\geq 0} A^{\bf S}_{+}(\alpha+q\beta) {\sf f}_{\alpha+q\beta} \big)
= \mathrm{exp}(-\IH)\, \mathrm{exp}\big(\sum_{q\geq 0} A^{\bf S}_{-}(\alpha+q\beta) {\sf f}_{\alpha+q\beta} \big)\,
\mathrm{exp}(\IH).
\ee
Now let us recall the following form of the BCH formula:
\begin{equation}
\begin{aligned}
\text{exp}(A) \text{exp}(B) \text{exp} (-A) & = \text{exp}( \sum_{n=0} \frac{1}{n!} (Ad(A))^n B )\\ & 
= \text{exp} ( B + [A,B] +\frac{1}{2} [A,[A,B]]+ \cdots)\\
\end{aligned}
\end{equation}
Using this formula in \eqref{eq:KSformD}, we obtain
\be\label{eq:KSformE}
\begin{aligned}
&\mathrm{exp} 
\big(\sum_{q\geq 0} A^{\bf S}_{+}(\alpha+q\beta) {\sf f}_{\alpha+q\beta} \big)
= \\
&\mathrm{exp}\big(\sum_{q\geq 0} A^{\bf S}_{-}(\alpha+q\beta) 
\sum_{l\geq 1}\mathop{\sum_{q_1,\ldots, q_l\geq 1}}_{}
{(-1)^l\over l!}\prod_{i=1}^l (-1)^{{\chi}(\alpha,q_i\beta)}{\chi}(\alpha,q_i\beta)
H(q_i\beta) {\sf f}_{\alpha+(q+q_1+\cdots + q_l)\beta}\big)\\
\end{aligned} 
\ee 
Finally, identifying the coefficients of a given Lie algebra generator 
${\sf f}_{\alpha+p\beta}$ we obtain the wallcrossing formula \eqref{eq:wallformulaAA}. 

\subsection{Comparison with Denef-Moore halo formula} 
Suppose $X$ is a genus zero curve such that the Higgs sheaf invariants $H(r,e)$ 
may be nontrivial.  
Employing the notation introduced in the previous subsection, consider the following generating functions
\[
Z_{\pm}(q,v) = \sum_{p\geq 0} A_{\delta_\pm}^{\bf S}(\alpha+p\beta) 
q^{pn(\beta)} v^{pr(\beta)},
\]
where $\alpha=(r(\alpha),e(\alpha))$, $\beta=(r(\beta),e(\beta))$ and $n(\beta) = 
e(\beta)-r(\beta)(g-1)$. 
Then the wallcrossing formula \eqref{eq:wallformulaAA} yields 
\[
\bal 
Z_+(q,v) = \sum_{l\geq 0} {(-1)^l\over l!} &
\bigg(\sum_{p\geq 0} A_{\delta_-}^{\bf S}(\alpha+p\beta)  q^{pn(\beta)} v^{pr(\beta)} \bigg) \\
& \bigg(\sum_{p\geq 0} (-1)^{\chi(\alpha, p\beta)} \chi(\alpha, p\beta) 
 H(p\beta) q^{pn(\beta)}v^{pr(\beta)} \bigg)^l.\\
\eal 
\]
Using the multicover formula \eqref{eq:higgsinvB}, 
this expression becomes 
\[ 
\bal 
Z_+(q,v) & = Z_-(q,v) \mathrm{exp}\bigg[ 
-\sum_{k,p\geq 0} {1\over k} (-1)^{kp\chi(\alpha, \beta)} \chi(\alpha, p\beta) 
 {\overline H}(p\beta) q^{kpn(\beta)}v^{kpr(\beta)}\bigg]\\
 & = Z_-(q,v) \mathrm{exp}\bigg[ \sum_{p\geq 0} \chi(\alpha, p\beta)  
 {\overline H}(p\beta)\mathrm{ln}\big( 
 1-(-1)^{p\chi(\alpha,\beta)} q^{pn(\beta)}\big) \bigg]\\ 
 & = Z_-(q,v)\prod_{p\geq 0} \big( 
 1-(-1)^{p\chi(\alpha,\beta)} q^{pn(\beta)}v^{pr(\beta)}\big)^{\chi(\alpha, p\beta) 
 {\overline H}(p\beta)}   \\
 \eal 
 \]
This formula in agreement with the halo formula \cite[Eqn. 6.17]{DM-split}.

  \appendix
  \section{Bell Polynomials} 
In this section we summarize some basic facts concerning Bell polynomials 
used in the proof of lemma (\ref{SFlemmD}) following \cite[Sect. 3.3-3.4]{advcomb}.

Let ${\sf P}_k(n)$ be the set of unordered length $k\geq 1$ partitions
of a positive integer $n\geq 1$. A partition $\lambda \in  {\sf P}_k(n)$ 
is determined by a sequence $(j_1,\ldots, j_{n-k+1})$ of non-negative 
integers satisfying 
\[
j_1+2j_2+\cdots =n, \qquad j_1+j_2 + \cdots = k.
\]
Then we write $\lambda =(1^{j_1}, 2^{j_2}, \ldots)$. 
For us, the  Bell polynomial $B_{n,k}(x_1,\ldots, x_{n-k+1})$ will be defined by the 
following formula 
\be\label{eq:BellA}
\begin{aligned}
& B_{n,k}(x_1,\ldots, x_{n-k+1}) =\\
& \mathop{\sum_{\lambda\in {\sf P}_k(n)}}_{}
{n!\over j_1!j_2!\cdots j_{n-k+1}!} {1\over (1!)^{j_1} (2!)^{j_2} \cdots 
((n-k+1)!)^{j_{n-k+1}}} x_1^{j_1}x_2^{j_2}\cdots x_{n-k+1}^{j_{n-k+1}}\\
\end{aligned}
\ee
The power series version of Fa\`a di Bruno's formula is the following identity
(see \cite[Thm.A, Sect. 3.4]{advcomb}. 
Let 
\[
f(x) = \sum_{n=1}^\infty {a_n\over n!} x^n \qquad 
 g(x) = \sum_{n=1}^\infty {b_n\over n!} x^n
\]
be formal power series with complex coefficients. 
Then 
\be\label{eq:FdB}
g(f(x)) = \sum_{n=1}^\infty {1\over n!} \sum_{k=1}^n b_kB_{n,k}(a_1,\ldots, a_{n-k+1}) 
x^n.
\ee
Now let 
\[
f(x) = e^{x}-1 = \sum_{n=1}^\infty {x^n\over n!} 
\qquad 
g(x) = -{x\over 1+x} = \sum_{n=1}^{\infty} (-1)^n x^n 
\]
Then
\[
g(f(x)) = -1-e^{-x} = \sum_{n=1}^\infty {(-1)^n\over n!} x^n,
\]
and equation \eqref{eq:FdB} yields 
\be\label{eq:Bellsum}
\sum_{k=1}^n (-1)^k k! \, B_{n,k}(1,\ldots,1) = (-1)^n.
\ee

\bibliography{newref.bib}
 \bibliographystyle{abbrv}
\end{document}